# Sup-norm convergence rate and sign concentration property of Lasso and Dantzig estimators

**Karim Lounici**

*Laboratoire de Statistiques CREST*
*3, avenue Pierre Larousse 92240 Malakoff France and*
*Laboratoire de Probabilités et Modèles Aléatoires (UMR CNRS 7599),*
*Université Paris 7, 2 pl. Jussieu, BP 7012, 75251 Paris Cedex 05, France.*
*e-mail:* `lounici@math.jussieu.fr`

**Abstract:** We derive the $l_\infty$ convergence rate simultaneously for Lasso and Dantzig estimators in a high-dimensional linear regression model under a mutual coherence assumption on the Gram matrix of the design and two different assumptions on the noise: Gaussian noise and general noise with finite variance. Then we prove that simultaneously the thresholded Lasso and Dantzig estimators with a proper choice of the threshold enjoy a sign concentration property provided that the non-zero components of the target vector are not too small.

## 1. Introduction

The Lasso is an $l_1$ penalized least squares estimator in linear regression models proposed by Tibshirani [17]. The Lasso enjoys two important properties. First, it is naturally sparse, i.e., it has a large number of zero components. Second, it is computationally feasible even for high-dimensional data (Efron et al. [8], Osborne et al. [16]) whereas classical procedures such as BIC are not feasible when the number of parameters becomes large. The first property raises the question of model selection consistency of Lasso, i.e., of identification of the subset of non-zero parameters. A closely related problem is sign consistency, i.e., identification of the non-zero parameters and their signs (cf. Bunea [2], Meinshausen and Bühlmann [13], Meinshausen and Yu [14], Wainwright [20], Zhao and Yu [22] and the references cited in these papers).

Zou [23] has proved estimation and variable selection results for the adaptive Lasso: a variant of Lasso where the weights on the different components in the $l_1$ penalty vary and are data dependent. We mention also work on the convergence of the Lasso estimator under the prediction loss: Bickel, Ritov and Tsybakov [1], Bunea, Tsybakov and Wegkamp [3], Greenshtein and Ritov [9], Koltchinskii [11; 12], Van der Geer [18; 19].

Knight and Fu [10] have proved the estimation consistency of the Lasso estimator in the case where the number of parameters is fixed and smaller than the sample size. The $l_2$ consistency of Lasso with convergence rate has been proved in Bickel, Ritov and Tsybakov [1], Meinshausen and Yu [14], Zhang and Huang [21]. These results trivially imply the $l_p$ consistency, with $2 \leqslant p \leqslant \infty$, however with a suboptimal rate (cf., e.g., Theorem 3 in [21]). Bickel, Ritov and Tsybakov [1] have proved that the Dantzig selector of Candes and Tao [6] shares a lot of common properties with the Lasso. In particular they have shown simultaneous $l_p$ consistency with rates of the Lasso and Dantzig estimators for $1 \leqslant p \leqslant 2$. To our knowledge, there is no result on the $l_\infty$ convergence rate and sign consistency of the Dantzig estimator.

The notion of $l_\infty$ and sign consistency should be properly defined when the number of parameters is larger than the sample size. We may have indeed an infinity of possible target vectors and solutions to the Lasso and Dantzig minimization problems. This difficulty is not discussed in [2; 13; 14; 20; 21] where either the target vector or the Lasso estimator or both are assumed to be unique. We show that under a sparsity scenario, it is possible to derive $l_\infty$ and sign consistency results even when the number of parameters is larger than the sample size. We refer to Theorem 6.3 and the Remark 1, p. 21, in [1] which suggest a way to clarify the difficulty mentioned above.

In this paper, we consider a high-dimensional linear regression model where the number of parameters can be much greater than the sample size. We show that under a mutual coherence assumption on the Gram matrix of the design, the target vector which has few non-zero components is unique. We do not assume the Lasso or Dantzig estimators to be unique. We establish the $l_\infty$ convergence rate of all the Lasso and Dantzig estimators simultaneously under two different assumptions on the noise. The rate that we get improves upon those obtained for the Lasso in the previous works. Then we show a sign concentration property of all the thresholded Lasso and Dantzig estimators simultaneously for a proper choice of the threshold if we assume that the non-zero components of the sparse target vector are large enough. Our condition on the size of the non-zero components of the target vector is less restrictive than in [20–22]. In addition, we prove analogous results for the Dantzig estimator, which to our knowledge was not done before.

The paper is organized as follows. In Section 2 we present the Gaussian linear regression model, the assumptions, the results and we compare them with the existing results in the literature. In Section 3 we consider a general noise with zero mean and finite variance and we show that the results remain essentially the same, up to a slight modification of the convergence rate. In Section 4 we provide the proofs of the results.

## 2. Model and Results

Consider the linear regression model

$$Y = X\theta^* + W, \tag{1}$$

where $X$ is an $n \times M$ deterministic matrix, $\theta^* \in \mathbb{R}^M$ and $W = (W_1, \ldots, W_n)^T$ is a zero-mean random vector such that $\mathbb{E}[W_i^2] \leqslant \sigma^2$, $1 \leqslant i \leqslant n$ for some $\sigma^2 > 0$. For any $\theta \in \mathbb{R}^M$, define $J(\theta) = \{j : \theta_j \neq 0\}$. Let $M(\theta) = |J(\theta)|$ be the cardinality of $J(\theta)$ and $\vec{\text{sign}}(\theta) = (\text{sign}(\theta_1), \ldots, \text{sign}(\theta_M))^T$ where

$$\text{sign}(t) = \begin{cases} 1 & \text{if } t > 0, \\ 0 & \text{if } t = 0, \\ 1 & \text{if } t < 0. \end{cases}$$

For any vector $\theta \in \mathbb{R}^M$ and any subset $J$ of $\{1, \ldots, M\}$, we denote by $\theta_J$ the vector in $\mathbb{R}^M$ which has the same coordinates as $\theta$ on $J$ and zero coordinates on the complement $J^c$ of $J$. For any integers $1 \leqslant d, p < \infty$ and $z = (z_1, \ldots, z_d) \in \mathbb{R}^d$, the $l_p$ norm of the vector $z$ is denoted by $|z|_p \triangleq \left( \sum_{j=1}^d |z_j|^p \right)^{1/p}$, and $|z|_\infty \triangleq \max_{1 \leqslant j \leqslant d} |z_j|$.

Note that the assumption of uniqueness of $\theta^*$ is not satisfied if $M > n$. In this case, if a vector $\theta^* = \theta^0$ satisfies (1), then there exists an affine space $\Theta^* = \{\theta^* : X\theta^* = X\theta^0\}$ of dimension $\geqslant M - n$ of vectors satisfying (1). So the question of sign consistency becomes problematic when $M > n$ because we can easily find two distinct vectors $\theta^1$ and $\theta^2$ satisfying (1) such that $\vec{\text{sign}}(\theta^1) \neq \vec{\text{sign}}(\theta^2)$. However we will show that under our assumption of sparsity $\theta^*$ is unique.

The Lasso and Dantzig estimators $\hat{\theta}^L, \hat{\theta}^D$ solve respectively the minimization problems

$$\min_{\theta \in \mathbb{R}^M} \frac{1}{n} |Y - X\theta|_2^2 + 2r|\theta|_1, \tag{2}$$

and

$$\min_{\theta \in \mathbb{R}^M} |\theta|_1 \text{ subject to } \left| \frac{1}{n} X^T (Y - X\theta) \right|_\infty \leqslant r, \tag{3}$$

where $r > 0$ is a constant. A convenient choice in our context will be $r = A\sigma\sqrt{(\log M)/n}$, for some $A > 0$. We denote respectively by $\hat{\Theta}^L$ and $\hat{\Theta}^D$ the set of solutions to the Lasso and Dantzig minimization problems (2) and (3).

The definition of the Lasso minimization problem we use here is not the same as the one in [17], where it is defined as

$$\min_{\theta \in \mathbb{R}^M} \frac{1}{n} |Y - X\theta|_2^2 \text{ subject to } |\theta|_1 \leqslant t,$$

for some $t > 0$. However these minimization problems are strongly related, cf. [5]. The Dantzig estimator was introduced and studied in [6]. Define $\Phi(\theta) = \frac{1}{n} |Y - X\theta|_2^2 + 2r|\theta|_1$. A necessary and sufficient condition for a vector $\theta$ to minimize $\Phi$ is that the zero vector in $\mathbb{R}^M$ belongs to the subdifferential of $\Phi$ at point $\theta$, i.e.,

$$\begin{cases} \frac{1}{n}(X^T(Y - X\theta))_j = \text{sign}(\theta_j)r & \text{if } \theta_j \neq 0, \\ \left| \frac{1}{n}(X^T(Y - X\theta))_j \right| \leqslant r & \text{if } \theta_j = 0. \end{cases}$$

Thus, any vector $\theta \in \hat{\Theta}^L$ satisfies the Dantzig constraint

$$\left|\frac{1}{n}X^T(Y - X\theta)\right|_\infty \leq r. \tag{4}$$

The Lasso estimator is unique if $M < n$, since in this case $\Phi(\theta)$ is strongly convex. However, for $M > n$ it is not necessarily unique. The uniqueness of Dantzig estimator is not granted either. From now on, we set $\hat{\Theta} = \hat{\Theta}^L$ or $\hat{\Theta}^D$ and $\hat{\theta}$ denotes an element of $\hat{\Theta}$.

Now we state the assumptions on our model. The first assumption concerns the noise variables.

**Assumption 1.** *The random variables $W_1, \ldots, W_n$ are i.i.d. $\mathcal{N}(0, \sigma^2)$.*

We also need assumptions on the Gram matrix

$$\Psi \triangleq \frac{1}{n}X^TX.$$

**Assumption 2.** *The elements $\Psi_{i,j}$ of the Gram matrix $\Psi$ satisfy*

$$\Psi_{j,j} = 1, \quad \forall 1 \leq j \leq M, \tag{5}$$

*and*

$$\max_{i \neq j} |\Psi_{i,j}| \leq \frac{1}{\alpha(1 + 2c_0)s}, \tag{6}$$

*for some integer $s \geq 1$ and some constant $\alpha > 1$, where $c_0 = 1$ if we consider the Dantzig estimator, and $c_0 = 3$ if we consider the Lasso estimator.*

The notion of mutual coherence was introduced in [7] where the authors required that $\max_{i \neq j} |\Psi_{i,j}|$ were sufficiently small. Assumption 2 is stated in a slightly weaker form in [1]-[4].

Consider two vectors $\theta^1$ and $\theta^2$ satisfying (1) such that $M(\theta^1) \leq s$ and $M(\theta^2) \leq s$. Denote $\theta = \theta^1 - \theta^2$ and $J = J(\theta^1) \cup J(\theta^2)$. We clearly have $X\theta = 0$ and $|J| \leq 2s$. Assume that $\theta \neq 0$. Under Assumption 2, similarly as we derive the inequality (11) in Section 4 below and using the fact that $|\theta|_1 \leq \sqrt{2s}|\theta|_2$, we get that

$$\frac{|X\theta|_2^2}{n|\theta|_2^2} > 0.$$

This contradicts the fact that $X\theta = 0$. Thus we have $\theta^1 = \theta^2$. We have proved that under Assumption 2 the vector $\theta^*$ satisfying (1) with $M(\theta^*) \leq s$ is unique.

Our first result concerns the $l_\infty$ rate of convergence of Lasso and Dantzig estimators.

**Theorem 1.** *Take $r = A\sigma\sqrt{(\log M)/n}$ and $A > 2\sqrt{2}$. Let Assumptions 1,2 be satisfied. If $M(\theta^*) \leq s$, then*

$$\mathbb{P}\left(\sup_{\hat{\theta} \in \hat{\Theta}} \left|\hat{\theta} - \theta^*\right|_\infty \leq c_2 r\right) \geq 1 - M^{1 - A^2/8},$$

*with $c_2 = \frac{3}{2}\left(1 + \frac{(1+c_0)^2}{(1+2c_0)(\alpha-1)}\right)$.*

Theorem 1 states that in high dimensions $M$ the set of estimators $\hat{\Theta}$ is necessarily well concentrated around the vector $\theta^*$. Similar phenomenon was already observed in [1], cf. Remark 1, page 21, for concentration in $l_p$ norms, $1 \leq p \leq 2$. Note that $c_2$ in Theorem 1 is an absolute constant. Using Theorem 1, we can easily prove the consistency of the Lasso and Dantzig estimators simultaneously when $n \to \infty$. We allow the quantities $s$, $M$, $\hat{\Theta}$, $\theta^*$ to vary with $n$. In particular, we assume that

$$M \to \infty \quad \text{and} \quad \lim_{n \to \infty} \frac{\log M}{n} = 0,$$

as $n \to \infty$, and that Assumptions 1,2 hold true for any $n$. Then we have

$$\sup_{\hat{\theta} \in \hat{\Theta}} \left| \hat{\theta} - \theta^* \right|_\infty \to 0 \tag{7}$$

in probability, as $n \to \infty$. The condition $(\log M)/n \to 0$ means that the number of parameters cannot grow arbitrarily fast when $n \to \infty$. We have the restriction $M = o(\exp(n))$, which is natural in this context.

A result on $l_\infty$ consistency of Lasso has been previously stated in Theorem 3 of [21], where $\hat{\theta}^L$ was assumed to be unique and under another assumption on the matrix $\Psi$. It is not directly related to our Assumption 2, but can be deduced from a restricted version of Assumption 2 where $\alpha$ is taken to be substantially larger than 1. The result in [21] is a trivial consequence of the $l_2$ consistency, and has therefore the rate $|\hat{\theta}^L - \theta^*|_\infty = O_\mathbb{P}(s^{1/2} r)$ which is slower than the correct rate given in Theorem 1. In fact, the rate in [21] depends on the unknown sparsity $s$ which is not the case in Theorem 1. Note also that Theorem 3 in [21] concerns the Lasso only, whereas our result covers simultaneously the Lasso and Dantzig estimators.

We now study the sign consistency. We make the following assumption.

**Assumption 3.** *There exists an absolute constant $c_1 > 0$ such that*

$$\rho \stackrel{\Delta}{=} \min_{j \in J(\theta^*)} |\theta_j^*| > c_1 r.$$

We will take $r = A\sigma \sqrt{(\log M)/n}$. We can find similar assumptions on $\rho$ in the work on sign consistency of the Lasso estimator mentioned above. More precisely, the lower bound on $\rho$ is of the order $s^{1/4} r^{1/2}$ in [14], $n^{-\delta/2}$ with $0 < \delta < 1$ in [20; 22], $\sqrt{(\log Mn)/n}$ in [2] and $\sqrt{s} r$ in [21]. Note that our assumption is the less restrictive.

We now introduce thresholded Lasso and Dantzig estimators. For any $\hat{\theta} \in \hat{\Theta}$ the associated thresholded estimator $\tilde{\theta} \in \mathbb{R}^M$ is defined by

$$\tilde{\theta}_j = \begin{cases} \hat{\theta}_j, & \text{if } |\hat{\theta}_j| > c_2 r, \\ 0 & \text{elsewhere.} \end{cases}$$

Denote by $\tilde{\Theta}$ the set of all such $\tilde{\theta}$. We have first the following non-asymptotic result that we call sign concentration property.

**Theorem 2.** *Take $r = A\sigma\sqrt{(\log M)/n}$ and $A > 2\sqrt{2}$. Let Assumptions 1-3 be satisfied. We assume furthermore that $c_1 > 2c_2$, where $c_2$ is defined in Theorem 1. Then*
$$\mathbb{P}\left(\vec{\operatorname{sign}}(\tilde{\theta}) = \vec{\operatorname{sign}}(\theta^*),\ \forall \tilde{\theta} \in \tilde{\Theta}\right) \geqslant 1 - M^{1-A^2/8}.$$

Theorem 2 guarantees that every vector $\tilde{\theta} \in \tilde{\Theta}$ and $\theta^*$ share the same signs with high probability. Letting $n$ and $M$ tend to $\infty$ we can deduce from Theorem 2 an asymptotic result under the following additional assumption.

**Assumption 4.** *We have $M \to \infty$ and $\lim_{n \to \infty} \frac{\log M}{n} = 0$, as $n \to \infty$.*

Then the following asymptotic result called sign consistency follows immediately from Theorem 2.

**Corollary 1.** *Let the assumptions of Theorem 2 hold for any $n$ large enough. Let Assumption 4 be satisfied. Then*
$$\mathbb{P}\left(\vec{\operatorname{sign}}(\tilde{\theta}) = \vec{\operatorname{sign}}(\theta^*),\ \ \forall \tilde{\theta} \in \tilde{\Theta}\right) \to 1,$$

*as $n \to \infty$.*

The sign consistency of Lasso was proved in [13; 22] with the Strong Irrepresentable Condition on the matrix $\Psi$ which is somewhat different from ours. Papers [13; 22] assume a lower bound on $\rho$ of the order $n^{-\delta/2}$ with $0 < \delta < 1$, whereas our Assumption 3 is less restrictive. Note also that these papers assume $\hat{\theta}^L$ to be unique. Wainwright [20] does not assume $\hat{\theta}^L$ to be unique and discusses sign consistency of Lasso under a mutual coherence assumption on the matrix $\Psi$ and the following condition on the lower bound: $\sqrt{(\log M)/n} = o(\rho)$ as $n \to \infty$, which is more restrictive than our Assumption 3. In particular Proposition 1 in [20] states that as $n \to \infty$, if the sequence of $\theta^*$ satisfies the above condition for all $n$ large enough, then
$$\mathbb{P}\left(\exists \hat{\theta}^L \in \hat{\Theta}^L \text{ s.t. } \vec{\operatorname{sign}}(\hat{\theta}^L) = \vec{\operatorname{sign}}(\theta^*)\right) \to 1.$$

This result does not guarantee sign consistency for all the estimators $\hat{\theta}^L \in \hat{\Theta}^L$ but only for some unspecified subsequence that is not necessarily the one chosen in practice. On the contrary, Corollary 1 guarantees that all the thresholded Lasso and Dantzig estimators and $\theta^*$ share the same sign vector asymptotically. It follows from this result that any solution selected by the minimization algorithm is covered and that the case $M > n$, where the set $\hat{\Theta}$ is not necessarily reduced to an unique estimator, can still be treated. We note also that the papers mentioned above treat the sign consistency for the Lasso only, whereas we prove it simultaneously for Lasso and Dantzig estimators. An improvement in the conditions that we get is probably due to the fact that we consider thresholded Lasso and Dantzig estimators. In addition note that not only the consistency results, but also the exact non-asymptotic bounds are provided by Theorems 1 and 2.

## 3. Convergence rate and sign consistency under a general noise

In the literature on Lasso and Dantzig estimators, the noise is usually assumed to be Gaussian [1; 6; 13; 20; 21] or admitting a finite exponential moment [2; 14]. The exception is the paper by Zhao and Yu [22] who proved the sign consistency of the Lasso when the noise admits a finite moment of order $2k$ where $k \geqslant 1$ is an integer. An interesting question is to determine whether the results of the previous section remain valid under less restrictive assumption on the noise. In this section, we only assume that the random variables $W_i, i = 1, \ldots, n$, are independent with zero mean and finite variance $\mathbb{E}[W_i^2] \leqslant \sigma^2$. We show that the results remain similar. We need the following assumption

**Assumption 5.** *The matrix $X$ is such that*

$$\frac{1}{n} \sum_{i=1}^{n} \max_{1 \leqslant j \leqslant M} |X_{i,j}|^2 \leqslant c',$$

*for a constant $c' > 0$.*

For example, if all $X_{i,j}$ are bounded in absolute value by a constant uniformly in $i, j$, then Assumption 4 is satisfied. The next theorem gives the $l_\infty$ rate of convergence of Lasso and Dantzig estimators under a mild noise assumption.

**Theorem 3.** *Assume that $W_i$ are independent random variables with $\mathbb{E}[W_i] = 0$, $\mathbb{E}[W_i^2] \leqslant \sigma^2$, $i = 1, \ldots, n$. Take $r = \sigma\sqrt{\frac{(\log M)^{1+\delta}}{n}}$, with $\delta > 0$. Let Assumptions 2,5 be satisfied. Then*

$$\mathbb{P}\left(\sup_{\hat{\theta} \in \hat{\Theta}} \left|\hat{\theta} - \theta^*\right|_\infty \leqslant c_2 r\right) \geqslant 1 - \frac{c}{(\log M)^\delta},$$

*where $c_2$ is defined in Theorem 1, and $c > 0$ is a constant depending only on $c'$.*

Therefore the $l_\infty$ convergence rate under the bounded second moment noise assumption is only slightly slower than the one obtained under the Gaussian noise assumption and the concentration phenomenon is less pronounced. If we assume that $\lim_{n \to \infty} (\log M)^{1+\delta}/n = 0$ and that Assumptions 2,3 and 5 hold true for any $n$ with $r = \sigma\sqrt{(\log M)^{1+\delta}/n}$, then the sign consistency of thresholded Lasso and Dantzig estimators follows from our Theorem 3 similarly as we have proved Theorem 2 and Corollary 1. Zhao and Yu [22] stated in their Theorem 3 a result on the sign consistency of Lasso under the finite variance assumption on the noise. They assumed $\hat{\theta}^L$ to be unique and the matrix $X$ to satisfy the condition $\max_{1 \leqslant i \leqslant n}(\sum_{j=1}^{M} X_{i,j}^2)/n \to 0$, as $n \to \infty$. This condition is rather strong. It does not hold if $M > n$ and all the $X_{i,j}$ are bounded in absolute value by a constant. In addition, [22] assumes that the dimension $M = O(n^\delta)$ with $0 < \delta < 1$, whereas we only need that $M = o(\exp(n^{1/(1+\delta)}))$ with $\delta > 0$. Note also that [22] proves the sign consistency for the Lasso only, whereas we prove it for thresholded Lasso and Dantzig estimators.

## 4. Proofs

We begin by stating and proving two preliminary lemmas. The first lemma originates from Lemma 1 of [3] and Lemma 2 of [1].

**Lemma 1.** *Let Assumption 1 and (5) of Assumption 2 be satisfied. Take $r = A\sigma\sqrt{(\log M)/n}$. Here $\hat{\Theta}$ denotes either $\hat{\Theta}^L$ or $\hat{\Theta}^D$. Then we have, on an event of probability at least $1 - M^{-A^2/8}$, that*

$$\sup_{\hat{\theta} \in \hat{\Theta}} \left|\Psi(\theta^* - \hat{\theta})\right|_\infty \leqslant \frac{3r}{2}, \qquad (8)$$

*and for all $\hat{\theta} \in \hat{\Theta}$,*

$$|\Delta_{J(\theta^*)^c}|_1 \leqslant c_0 |\Delta_{J(\theta^*)}|_1, \qquad (9)$$

*where $\Delta = \hat{\theta} - \theta^*$, $c_0 = 1$ for the Dantzig estimator and $c_0 = 3$ for the Lasso.*

*Proof.* Define the random variables $Z_j = n^{-1}\sum_{i=1}^n X_{i,j}W_i$, $1 \leqslant j \leqslant M$. Using (5) we get that $Z_j \sim \mathcal{N}(0, \sigma^2/n)$, $1 \leqslant j \leqslant M$. Define the event

$$\mathcal{A} = \bigcap_{j=1}^M \{|Z_j| \leqslant r/2\}.$$

Standard inequalities on the tail of Gaussian variables yield

$$\begin{aligned} P(\mathcal{A}^c) &\leqslant MP(|Z_1| \geqslant r/2), \\ &\leqslant M\exp\left(-\frac{n}{2\sigma^2}\left(\frac{r}{2}\right)^2\right) \\ &\leqslant M^{1-\frac{A^2}{8}}. \end{aligned}$$

On the event $\mathcal{A}$, we have

$$\left|\frac{1}{n}X^TW\right|_\infty \leqslant \frac{r}{2}. \qquad (10)$$

Any vector $\hat{\theta}$ in $\hat{\Theta}^L$ or $\hat{\Theta}^D$ satisfies the Dantzig constraint (4). Thus we have on $\mathcal{A}$ that

$$\sup_{\hat{\theta} \in \hat{\Theta}} \left|\Psi(\theta^* - \hat{\theta})\right|_\infty \leqslant \frac{3r}{2}.$$

Now we prove the second inequality. For any $\hat{\theta}^D \in \hat{\Theta}^D$, we have by definition that $|\hat{\theta}^D|_1 \leqslant |\theta^*|_1$, thus

$$|\Delta_{J(\theta^*)^c}|_1 = \sum_{j \in J(\theta^*)^c} |\hat{\theta}_j^D| \leqslant \sum_{j \in J(\theta^*)} |\theta_j^*| - |\hat{\theta}_j^D| \leqslant |\Delta_{J(\theta^*)}|_1.$$

Consider now the Lasso estimators. By definition, we have for any $\hat{\theta}^L \in \hat{\Theta}^L$

$$\frac{1}{n}|Y - X\hat{\theta}^L|_2^2 + 2r|\hat{\theta}^L|_1 \leqslant \frac{1}{n}|W|_2^2 + 2r|\theta^*|_1.$$

Developing the left hand side on the above inequality, we get

$$2r|\hat{\theta}^L|_1 \leqslant 2r|\theta^*|_1 + \frac{2}{n}(\hat{\theta}^L - \theta^*)^T X^T W.$$

On the event $\mathcal{A}$, we have for any $\hat{\theta}^L \in \hat{\Theta}^L$

$$2|\hat{\theta}^L|_1 \leqslant 2|\theta^*|_1 + |\hat{\theta}^L - \theta^*|_1,$$

Adding $|\hat{\theta}^L - \theta^*|_1$ on both side, we get

$$\begin{aligned}|\hat{\theta}^L - \theta^*|_1 + 2|\hat{\theta}^L|_1 &\leqslant 2|\theta^*|_1 + 2|\hat{\theta}^L - \theta^*|_1 \\ |\hat{\theta}^L - \theta^*|_1 &\leqslant 2(|\hat{\theta}^L - \theta^*|_1 + |\theta^*|_1 - |\hat{\theta}^L|_1),\end{aligned}$$

Now we remark that if $j \in J(\theta^*)^c$, then we have $|\hat{\theta}^L_j - \theta^*_j| + |\theta^*_j| - |\hat{\theta}^L_j| = 0$. Thus we have on the event $\mathcal{A}$ that

$$\begin{aligned}|\Delta_{J(\theta^*)^c}|_1 - |\Delta_{J(\theta^*)}|_1 &\leqslant |\Delta|_1 \leqslant 2|\Delta_{J(\theta^*)}|_1 \\ |\Delta_{J(\theta^*)^c}|_1 &\leqslant 3|\Delta_{J(\theta^*)}|_1,\end{aligned}$$

for any $\hat{\theta}^L \in \hat{\Theta}^L$. □

**Lemma 2.** *Let Assumption 2 be satisfied. Then*

$$\kappa(s, c_0) \triangleq \min_{J \subset \{1,\cdots,M\}, |J| \leqslant s} \min_{\lambda \neq 0: |\lambda_{J^c}|_1 \leqslant c_0 |\lambda_J|_1} \frac{|X\lambda|_2}{\sqrt{n}|\lambda_J|_2} \geqslant \sqrt{1 - \frac{1}{\alpha}} > 0.$$

*Proof.* For any subset $J$ of $\{1, \ldots, M\}$ such that $|J| \leqslant s$ and $\lambda \in \mathbb{R}^M$ such that $|\lambda_{J^c}|_1 \leqslant c_0|\lambda_J|_1$, we have

$$\begin{aligned}\frac{|X\lambda_J|_2^2}{n|\lambda_J|_2^2} &= 1 + \frac{\lambda_J^T(\Psi - I_M)\lambda_J}{|\lambda_J|_2^2} \\ &\geqslant 1 - \frac{1}{\alpha(1+2c_0)s} \sum_{i,j=1}^{M} \frac{|\lambda_J^{(i)}||\lambda_J^{(j)}|}{|\lambda_J|_2^2} \\ &\geqslant 1 - \frac{1}{\alpha(1+2c_0)s} \frac{|\lambda_J|_1^2}{|\lambda_J|_2^2}, \quad (11)\end{aligned}$$

where we have used Assumption 2 in the second line, $I_M$ denotes the $M \times M$ identity matrix and $\lambda_J = (\lambda_J^{(1)}, \ldots, \lambda_J^{(M)})$ denotes the components of the vector $\lambda_J$. This yields

$$\begin{aligned}\frac{|X\lambda|_2^2}{n|\lambda_J|_2^2} &\geqslant \frac{|X\lambda_J|_2^2}{n|\lambda_J|_2^2} + 2\frac{\lambda_J^T X^T X \lambda_{J^c}}{n|\lambda_J|_2^2} \\ &\geqslant 1 - \frac{1}{\alpha s(1+2c_0)} \frac{|\lambda_J|_1^2}{|\lambda_J|_2^2} - \frac{2}{\alpha s(1+2c_0)} \frac{|\lambda_J|_1 |\lambda_{J^c}|_1}{|\lambda_J|_2^2} \\ &\geqslant 1 - \frac{1}{\alpha s} \frac{|\lambda_J|_1^2}{|\lambda_J|_2^2} \\ &\geqslant 1 - \frac{1}{\alpha} > 0.\end{aligned}$$

We have used Assumption 2 in the second line, the inequality $|\lambda_{J^c}|_1 \leqslant c_0|\lambda_J|_1$ in the third line and the fact that $|\lambda_J|_1 \leqslant \sqrt{|J|}|\lambda_J|_2 \leqslant \sqrt{s}|\lambda_J|_2$ in the last line. $\square$

**Proof of Theorem 1.** For all $1 \leqslant j \leqslant M$, $\hat{\theta} \in \hat{\Theta}$ we have

$$(\Psi(\theta^* - \hat{\theta}))_j = (\theta_j^* - \hat{\theta}_j) + \sum_{i=1, i \neq j}^{M} \Psi_{i,j}(\theta_i^* - \hat{\theta}_i).$$

Assumption 2 yields

$$|(\Psi(\theta^* - \hat{\theta}))_j - (\theta_j^* - \hat{\theta}_j)| \leqslant \frac{1}{\alpha(1+2c_0)s} \sum_{i=1, i \neq j}^{M} |\theta_i^* - \hat{\theta}_i|, \forall j.$$

Thus we have

$$|\theta^* - \hat{\theta}|_\infty \leqslant \left|\Psi(\theta^* - \hat{\theta})\right|_\infty + \frac{1}{\alpha(1+2c_0)s}|\theta^* - \hat{\theta}|_1. \tag{12}$$

Set $\Delta = \hat{\theta} - \theta^*$. Lemma 1 yields that on an event $\mathcal{A}$ of probability at least $1 - M^{1-A^2/8}$ we have for any $\hat{\theta} \in \hat{\Theta}$

$$|\Psi\Delta|_\infty \leqslant \frac{3r}{2}, \tag{13}$$

and

$$|\Delta|_1 = |\Delta_{J(\theta^*)^c}|_1 + |\Delta_{J(\theta^*)}|_1 \leqslant (1+c_0)|\Delta_{J(\theta^*)}|_1 \leqslant (1+c_0)\sqrt{s}|\Delta_{J(\theta^*)}|_2.$$

Thus we have, on the same event $\mathcal{A}$,

$$\begin{aligned}\frac{1}{n}|X\Delta|_2^2 &= \Delta^T\Psi\Delta \\ &\leqslant |\Psi\Delta|_\infty|\Delta|_1 \\ &\leqslant \frac{3r}{2}(1+c_0)\sqrt{s}|\Delta_{J(\theta^*)}|_2,\end{aligned} \tag{14}$$

for any $\hat{\theta} \in \hat{\Theta}$. Lemma 2 yields

$$\frac{1}{n}|X\Delta|_2^2 \geqslant \left(1 - \frac{1}{\alpha}\right)|\Delta_{J(\theta^*)}|_2^2, \tag{15}$$

for any $\hat{\theta} \in \hat{\Theta}$. Combining (14) and (15), we obtain that

$$|\Delta|_1 \leqslant \frac{3}{2}r(1+c_0)^2\frac{\alpha}{\alpha-1}s, \tag{16}$$

for any $\hat{\theta} \in \hat{\Theta}$. Combining (12), (13) and (16) we obtain that

$$\sup_{\hat{\theta} \in \hat{\Theta}} |\hat{\theta} - \theta^*|_\infty \leqslant \frac{3}{2}\left(1 + \frac{(1+c_0)^2}{(1+2c_0)(\alpha-1)}\right)r. \square$$

**Proof of Theorem 2.** Theorem 1 yields $\sup_{\hat{\theta}\in\hat{\Theta}}|\hat{\theta}-\theta^*|_\infty \leqslant c_2 r$ on an event $\mathcal{A}$ of probability at least $1-M^{1-A^2/8}$. Take $\hat{\theta}\in\hat{\Theta}$. For $j\in J(\theta^*)^c$, we have $\theta_j^*=0$, and $|\hat{\theta}_j|\leqslant c_2 r$ on $\mathcal{A}$. For $j\in J(\theta^*)$, we have by Assumption 3 that $|\theta_j^*|\geqslant c_1 r$ and $|\theta_j^*|-|\hat{\theta}_j|\leqslant |\theta_j^*-\hat{\theta}_j|\leqslant c_2 r$ on $\mathcal{A}$. Since we assume that $c_1 > 2c_2$, we have on $\mathcal{A}$ that $|\hat{\theta}_j|\geqslant (c_1-c_2)r > c_2 r$. Thus on the event $\mathcal{A}$ we have: $j\in J(\theta^*) \Leftrightarrow |\hat{\theta}_j| > c_2 r$. This yields $\text{sign}(\tilde{\theta}_j) = \text{sign}(\hat{\theta}_j) = \text{sign}(\theta_j^*)$ if $j\in J(\theta^*)$ on the event $\mathcal{A}$. If $j\notin J(\theta^*)$, $\text{sign}(\theta_j^*) = 0$ and $\tilde{\theta}_j = 0$ on $\mathcal{A}$, so that $\text{sign}(\tilde{\theta}_j) = 0$. The same reasoning holds true simultaneously for all $\hat{\theta}\in\hat{\Theta}$ on the event $\mathcal{A}$. Thus we get the result. $\square$

**Proof of Theorem 3.** The proof of Theorem 3 is similar to the one of Theorem 1 up to a modification of the bound on $P(\mathcal{A}^c)$ in Lemma 1. Recall that $Z_j = n^{-1}\sum_{i=1}^n X_{i,j}W_i$, $1\leqslant j\leqslant M$ and the event $\mathcal{A}$ is defined by

$$\mathcal{A} = \bigcap_{j=1}^M \{|Z_j|\leqslant r/2\} = \{\max_{1\leqslant j\leqslant M}|Z_j|\leqslant r/2\}.$$

The Markov inequality yields that

$$P(\mathcal{A}^c) \leqslant \frac{4\mathbb{E}[\max_{1\leqslant j\leqslant M} Z_j^2]}{r^2}.$$

Then we use Lemma 3 given below with $p=\infty$ and the random vectors $Y_i = (X_{i,1}W_i/n,\ldots,X_{i,M}W_i/n)\in\mathbb{R}^M$, $i=1,\ldots,n$. We get that

$$P(\mathcal{A}^c) \leqslant \tilde{c}\frac{\log M}{r^2}\sigma^2 \sum_{i=1}^n \max_{1\leqslant j\leqslant M}\frac{X_{i,j}^2}{n^2},$$

where $\tilde{c} > 0$ is an absolute constant. Taking $r = \sigma\sqrt{(\log M)^{1+\delta}/n}$ and using Assumption 5 yields that

$$P(\mathcal{A}^c) \leqslant \frac{c}{(\log M)^\delta},$$

where $c > 0$ is an absolute constant. $\square$

The following result is Lemma 5.2.2, page 188 of [15].

**Lemma 3.** *Let $Y_1,\ldots,Y_n \in \mathbb{R}^M$ be independent random vectors with zero means and finite variance, and let $M\geqslant 3$. Then for every $p\in[2,\infty]$, we have*

$$\mathbb{E}\left[|\sum_{i=1}^n Y_i|_p^2\right] \leqslant \tilde{c}\min[p,\log M]\sum_{i=1}^n \mathbb{E}\left[|Y_i|_p^2\right],$$

*where $\tilde{c} > 0$ is an absolute constant.*

**Acknowledgements**

I wish to thank my advisor, Alexandre Tsybakov, for insightful comments and the time he kindly devoted to me.